\documentclass[12pt,reqno]{article}
\usepackage{amsfonts}
\usepackage{amsthm}
\usepackage[mathscr]{eucal}
\usepackage{amsmath}
\usepackage{amssymb}
\usepackage{amscd}
\swapnumbers

\theoremstyle{plain}
\newtheorem{lemma}{Lemma}[section]
\newtheorem{theorem}[lemma]{Theorem}
\newtheorem{proposition}[lemma]{Proposition}
\newtheorem{corollary}[lemma]{Corollary}
\theoremstyle{definition}
\newtheorem{example}[lemma]{Example}

\newtheorem{remark}[lemma]{Remark}

\def\d{$\displaystyle}
\def\be{\begin{equation}}
\def\heavycdot{\raisebox{2pt}{\tiny\kern.5em$\bullet$}}

\numberwithin{equation}{section}
\def\begeq{\stepcounter{lemma}\begin{equation}}

\date{}

\begin{document}

\title{ON CYCLIC HOMOLOGY OF $A_\infty$-ALGEBRAS}
\author{Masoud Khalkhali\thanks{\; The author is grateful to J. Williams
for typing the manuscript.  The author is supported by NSERC of
Canada.}}

\maketitle

%%%%%%%%%%%%%% abstract

%\begin{abstract}
%%%%%%% put the text of the abstract here
%\end{abstract}

The category of $A_\infty$-algebras extend the category of differential
graded $(DG)$ algebras.  The main result of the present paper asserts
that the periodic cyclic homology $HP_\bullet(A,m)$ of an $A_\infty$-algebra
$(A,m)$ is equal to ordinary periodic cyclic homology $HP_\bullet(H_0(A))$ of
the homology of $(A,m_1)$ in degree zero.  This result extends a well
known result of T. Goodwillie [Go] for the periodic cyclic homology of
$DG$ algebras.  We notice, however, that while the study of, at least
certain aspects of, cyclic homology of $DG$ algebras can be ``reduced''
to the study of cyclic homology of algebras by carefully adapting the
algebra case to cyclic objects in the category of chain complexes as in
[Go], no such technique is available for $A_\infty$-algebras.  In fact,
although this is perhaps possible, and certainly interesting to try, we
are not aware of any generalization of Connes' cyclic category [C3] to a
hypothetical ``$\infty$-cyclic category'', so that an $A_\infty$-algebra
yields a ``$\infty$-cyclic'' object.

Hochschild and cyclic homology of $A_\infty$-algebras were first defined
by Getzler and Jones in [GJ], where a (b,B) type bicomplex was defined
for an $A_\infty$-algebra.  Our approach is however different and is
based on ideas of Cuntz and Quillen [CQ$_1$,CQ$_2$,Q] and especially their
$X$-complex approach to cyclic homology.  In fact this is crucial for
the Cartan homotopy formula that we need.  As is evident in the present
paper, this approach allows a unified treatment of cyclic homology type
theories for various algebraic structures.  A key element in the proof of
our main theorem is a Cartan homotopy formula for quasifree $DG$ (co)
algebras from [Kh], which in turn is a generalization of the Cartan
homotopy formula of Cuntz and Quillen [CQ$_1$].

The concepts of $A_\infty$-spaces and $A_\infty$-algebras are due to
Stasheff in [S], where it is shown that a topological space has the
homotopy type of a (based) loop space if and only if it is an
$A_\infty$-space.  The Moore model of a loop space has the extra property
that it has a strictly associative product and hence its singular chains
is a $DG$ algebra.  This is used in [Go] to identify the
$SO(2)$-equivariant homology of the free loop space with the cyclic
homology of based loop space.  It follows from Proposition 2.3 that this
is indeed true for any model of the loop space.

Apart from applications in [GJ] to the $SO(2)$-equivariant version of
Chen's iterated integrals, the new surge of interest in
$A_\infty$-algebras goes back in part to a paper of M. Kontsevich [K],
where it is shown that cyclic cohomology classes in degree zero of any
$A_\infty$-algebra can be used to construct homology classes on the
moduli spaces of algebraic curves (see also [PS] where this construction
is further explained).

\section {The bar construction for $A_\infty$-algebras}
The relevance of bar construction to the study of cyclic homology of
associative algebras is explained in [Q].  Since a similar construction
is available for $A_\infty$-algebras [S,GJ], it is natural to expect
that it would play a similar role for Hochschild and cyclic homology of
$A_\infty$-algebras.  In general, the bar construction makes the
definitions and the identities to be satisfied by various operators in an
$A_\infty$-algebra setting completely natural.

Let us fix our notations and recall some elementary concepts from
``graded mathematics''.  Throughout this paper we work over a field $K$
of characteristic zero and all homomorphsims and tensor products are
over $K$.  If \d A = \bigoplus_{i\in\mathbb{Z}}A_i$ is a
$\mathbb{Z}$-graded vector space, we deonte by $A[n]$ the shifted graded
vector space defined by $A[n]_i = A_{i-n}$, $i\in\mathbb{Z}$.  We have
$A[n][m] = A[n+m]$.  The degree of a homogeneous element $a$ in a graded
space is denoted by $|a|$.  In case of several gradings, we will make
the required grading explicit.

A homogeneous map $f:A\to B$ between graded spaces has degree
$k\in\mathbb{Z}$, if $f(A_i)\subset B_{i+k}$ $\forall\;
i\in\mathbb{Z}$.  We have
\[ Hom^k(A,B) = \prod_{i\in\mathbb{Z}} Hom(A_i,B_{i+k})\; . \]
The tensor product of graded spaces $A$ and $B$ is defined by \d (A\otimes
B)_n = \bigoplus_{i+j=n}A_i\otimes B_j$.  If $f:A\to A^\prime$, $g:B\to
B^\prime$ are graded maps, then their graded tensor product $f\otimes
g:A\otimes B\to A^\prime\otimes B^\prime$ is defined by
\[ (f\otimes g)(a\otimes b) = (-1)^{|a|\; |g|}f(a)\otimes g(b)\; . \]
Note that $|f\otimes g| = |f|+|g|$.  The graded twist $\sigma :A\otimes B\to
B\otimes A$ is defined by $\sigma (a\otimes b) = (-1)^{|a|\; |b|}b\otimes
a$.  We will denote a tensor $a_1\otimes\dots\otimes a_n$ in $A^{\otimes
n}$ by $(a_1,\dots ,a_n)$.  Note that with our conventions, \d |(a_1,\dots
,a_n)| = \sum^n_{i=1}|a_i|$ and hence the r-th components of $A^{\otimes n}$
is given by \d (A^{\otimes n})_r = \{ (a_1,\dots , a_s)\; |\;\Sigma |a_i| =
r\} = \bigoplus_{s\geq0}\bigoplus_{i_1+\dots +i_s=r} A_{i_1}\otimes
A_{i_2}\otimes\dots\otimes A_{i_s}$.

Let $T^cA$ denote the cofree coaugmented graded coalgebra
generated by the positively graded vector space \d A =
\bigoplus_{i\geq 0}A_i$.  We have $(T^cA)_n=A^{\otimes n}$.  Explicitly,
\begin{align*}
& (T^cA)_0 = K\oplus A_0\oplus A^{\otimes 2}_0\oplus A^{\otimes
3}_0\oplus\dots\\
& (T^cA)_1 = A_1\oplus A_0\otimes A_1\oplus A_1\otimes A_0\oplus\dots\\
& (T^cA)_2 = A_2\oplus A_1\otimes A_1\oplus A_0\otimes A_1\otimes
A_0\oplus\dots\\
& \quad\vdots
\end{align*}
The coproduct $\Delta :T^cA\to T^cA\otimes T^cA$, which is a degree zero
map, is defined by $\Delta (1)=1\otimes 1$ and
\[ \Delta (a_1,\dots , a_n) = 1\otimes (a_1,\dots
,a_n)+\sum^n_{i=1}(a_1,\dots ,a_i)\otimes (a_{i+1},\dots ,a_n) +
(a_1,\dots ,a_n)\otimes 1\; . \]

Let $C$ be a graded counital coalgebra.  By a graded $C$-bicomodule we
mean a graded vector space $M$ such that the left and right coactions
$\Delta_\ell :C\to C\otimes M$ and $\Delta_r:C\to M\otimes C$ are of degree
zero.  We further assume that our bicomuldes are counitary.  A {\em
graded coderivation} of degree $k\in\mathbb{Z}$ is a degree $k$ map
$\delta :M\to C$ such that
\[ (1\otimes\delta )\Delta_\ell+(\delta\otimes 1)\Delta_r=\Delta\delta\; .
\]

A {\em universal coderivation} of degree $k$ consists of the following
data:  A $C$-bicomodule $\Omega^1_kC$ and a graded coderivation of
degree $k$, $d_k:\Omega^1_kC\to C$ which is universal.  That is, for any
degree $k$ coderivation $\delta :M\to C$, there exists a unique degree
zero $C$-bicomodule map $m:M\to\Omega^1_kC$ such that $\delta =
d_k\circ m$.  We refer to $\Omega^1_kC$ as the comodule of universal
codifferentials over $C$.

It is not difficult to see that universal graded coderivations of any
degree exist.  One can simply define $\Omega^1_kC =
coker\{\Delta_k:C[-k]\to C\otimes C\}$ with its $C$-bicomodule structure
induced from $C\otimes C$.  Here $\Delta_k$ is the composition $C[-k]\to
C\overset{\Delta}{\longrightarrow}C\otimes C$ and $d_k = \eta\otimes
1-1\otimes\eta$, where $\eta :C\to K$ is the counit.

In the cofree case $C=T^cA$, this universal coderivation can be
identified as follows.  Let
\[ \Omega^1_kT^cA = T^cA\otimes A[-k]\otimes T^cA \]
be the free bicomodule over $T^cA$ generated by $A[-k]$.  Define
$d_k:\Omega^1_kT^cA\to T^cA$ by
\[ d_k(\alpha\otimes a\otimes\beta) = (-1)^{|\alpha |}(\alpha ,a,\beta
), \]
where $\alpha ,\beta\in T^cA$ and $a\in A[-k]$.  It is not difficult to
see that $d_k$ is a degree $k$ coderivation and moreover it is
universal.

If in the adjunction formula
\[ Coder^k(M,T^cA)\simeq Hom(M,\Omega^1_kT^cA)\; , \]
we take $M=C$, with its natural bocomodule structure, we obtain, for
$k\in\mathbb{Z}$,
\[ Coder^k(T^cA,T^cA)\simeq Hom^k(T^cA,\Omega^1_kT^cA)\; . \]
However, $\Omega^1_kT^cA = T^cA\otimes A[-k]\otimes T^cA$, is a cofree
$T^cA$-bicomodule, and hence
\[ Hom^k(T^cA,\Omega^1_kT^cA)\simeq Hom^0(T^cA,A[-k])\simeq
Hom^k(T^cA,A)\; , \]
so, that we have an isomorphism of graded vectors spaces
\begin{equation}
Coder(T^cA,T^cA)\simeq Hom(T^cA,A)\; .\tag{1}
\end{equation}

This isomorphism works as follows.  A graded coderivation $b^\prime
:T^cA\to T^cA$ of degree $k$ defines a degree $k$ map
$m_{b^\prime}:T^cA\to A$ as the composition
\[ T^cA\overset{b^\prime}{\longrightarrow}T^cA\longrightarrow A\; , \]
where the last map is the projection.
Conversely, given a degree $k$ map $m:T^cA\to A$, define a degree $k$
coderivation $b^\prime_m:T^cA\to T^cA$ by
\[ b^\prime_m = d_k(1\otimes m\otimes 1)\; (1\otimes\Delta )\Delta\; .
\]
Chasing the formulas for $\Delta$ and $d_k$, we obtain the formula for
$b^\prime_m$.  We have $b^\prime_m(1) = 0$ and
\begin{eqnarray*}
b^\prime_m(a_1,\dots , a_n) & = & \sum^n_{i=1}b^\prime_{m_i}(a_1,\dots
,a_n)\\
& = &
\sum^n_{i=1}\left(\sum^{n-i}_{j=1}(-1)^{\varepsilon_{ij}}(a_1,\dots
,m_i(a_j,\dots ,a_{j+i}),\dots ,a_n)\right)\; .
\end{eqnarray*}

Finally, note that $T^cA$ has also a universal property with respect to
morphisms of graded coalgebras.  Namely, for any graded coaugmented
counital coalgebra $C$, we have an isomorphism of graded vector spaces
\begin{equation*}
 Hom^{GC}(C,T^cA) \simeq Hom(C,A)\; ,\tag{2}
\end{equation*}
where on the left hand side $Hom^{GC}$ means coaugmented graded
coalgebra morphisms.

Under this isomorphism, a graded linear map $f:C\to A$ of degree $k$
defines a coalgebra map $\hat{f}:C\to T^cA$, where its degree $n$ component
is given by $\hat{f}_n = f^{\otimes n}\circ\Delta^{(n)}_c$.  Here
$\Delta^{(n)}$ denotes the $n^{th}$ iteration of the comultiplication
$\Delta$.

Next we turn to the cofree coalgebra $C = T^cA[1]$, generated by the
suspension $A[1]$, where, as before $A$ is positively graded.  The reason
for this is twofold.  First of all, similar to algebra case [Q], various
cyclic bicomplexes for an $A_\infty$-algebra are obtained from
$T^cA[1]$.  Secondly, and independently, the all important Gerstenhaber
product that we are going to define preserves only the total grading
(= length + degree), which is the same as the grading in $A[1]$.  We have
\begin{align*}
& C_0 = K\\
& C_1 = A_0,\\
& C_2 = A_1\oplus A_0\otimes A_0,\\
& C_3 = A_2\oplus A_0\otimes A_1\oplus A_1\otimes A_0,\\
& \quad\vdots
\end{align*}
In general
\[ (T^cA[1])_n =\bigoplus^n_{r=0}\;\bigoplus_{i_1+\dots +i_r+r=n}
A_{i_1}\otimes\dots\otimes A_{i_r}\; . \]

A linear map $m:T^cA[1]\to A[1]$ of degree $k$ such that $m(1) = 0$ is
given by a sequence of linear maps
\[ m_n:A^{\otimes n}\to A\quad n=1,2,\dots \]
such that $|m_n| = n-1+k$.  We refer to $k$ as the suspended degree of
$m$.

The {\em Gerstenhaber product} (see [G] for the non-graded version) is a
degree zero non-associative product
\[ Hom(T^cA[1],A[1])\otimes Hom(T^cA[1],A[1])\to Hom(T^cA[1],A[1])\; ,
\]
defined as follows.  First, for $m:A^{\otimes k+1}\to A$, $m^\prime
:A^{\otimes\ell +1}\to A$, with suspended degrees $|m|$ and $|m^\prime
|$, define
\[ m\circ m^\prime :A^{\otimes k+\ell +1}\to A\; , \]
by
\[ (m\circ m^\prime )(a_1,\dots ,a_{k+\ell +1}) =
\sum^k_{i=1}(-1)^{\varepsilon_i}m(a_1,\dots ,a_{i-1},m^\prime (a_i,\dots
,a_{i+\ell}),\dots ,a_{k+\ell +1})\; , \]
where \d \varepsilon_i = |m^\prime
|\left(\sum^{i-1}_{j=1}|a_j|+i-1\right)$.  Note that we have the equality
of suspended degrees $|m\circ m^\prime |=|m|+|m^\prime |$.  We extend
this product to arbitrary cochains \d m = \sum^\infty_{i=1}m_i$ and \d
m^\prime = \sum^\infty_{i=1}m^\prime_i\in Hom(T^cA[1],A[1])$ by
\[ m\circ m^\prime = \sum^\infty_{n=2}\left(\sum_{i+j=n}m_i\circ
m_j\right)\; . \]

In [G], where the non-graded case is treated, Gerstenhaber proves that
$[m,m^\prime ]: = m\circ m^\prime - (-1)^{|m|\; |m^\prime |}m^\prime\circ m$
is a graded Lie bracket on $Hom(T^cA[1],A[1])$.  This is a surprise,
given the fact that $\circ$ is not associative.  This led to investigating
the full structure of higher homotopies that is hidden here [GV] (see
also section 3 of the present article).  The Lie algebra structure
itself, however, can be understood using isomorphism (1).  In fact, since
the bracket of two (graded) (co)derivation is again a
(graded)(co)derivation, the Lie strcture on the left hand side of (1) is
obvious and hence suffices to show that (1) preserves the brackets.

Now a map $m:T^cA[1]\to A[1]$ of degree $-1$ defines a degree $-1$
coderivation
\[ b^\prime_m:T^cA[1]\to T^cA[1]\; . \]
We have $b^{\prime^2}_m = \frac{1}{2}[b^\prime_m,b^\prime_m] =
\frac{1}{2}b^\prime_{[m,m]} = b^\prime_{m\circ m}$, so that $b^\prime_m$
is a codifferential iff $m\circ m = 0$.  Writing \d m =
\sum^\infty_{i=1}m_i$, we have $m\circ m=0$ iff
\[ \sum_{i+j=n}m_i\circ m_j = 0\qquad n=2,3,\dots\; . \]

An {\em $A_\infty$-algebra structure} (also called {\em strongly
homotopy associative algebra structure}) on a graded vector space \d A =
\bigoplus^\infty_{i=0} A_i$, is a degree $-1$ map $m:T^cA[1]\to A[1]$
such that $m\circ m = 0$.  Equivalently, it is defined by a coderivation
$b^\prime_m:T^cA[1]\to T^cA[1]$ of degree $-1$, such that
$b^{\prime^2}_m = 0$.  This concept is due to Stasheff [S].

Note that if $(A,m)$ is an $A_\infty$-algebra, the homotopy associative
product induces a strictly associative product on \d H_\bullet(A,m): =
H_\bullet(A,m_1) = \bigoplus_{i\geq 0}H_i(A,m_1)$ and turns it into an
associative graded algebra.  In particular $H_0(A,m_1)$ is an associative
algebra.

A morphism $(A,m)\to (B,m)$ of $A_\infty$-algebras is, by definition
[GJ], a morphism
\[ (T^cA[1],b^\prime_m)\to (T^cB(1],b^\prime_m) \]
of the corresponding coaugmented DG coalgebras.  The universal
property (2) has an obvious extension from graded vector spaces to complexes.
Using this, we see that there is a 1-1 correspondence between morphisms
of $A_\infty$-algebras $(A,m)\to (B,m)$ and morphisms of complexes
\[ f:(T^cA[1],b^\prime_m)\to (B[1],m_1)\; . \]

This notion of morphism between $A_\infty$-algebras may seem too
general, but, as we will see, all homological invariants that we
construct are in fact functorial with respect to these morphisms.  A
{\em strict morphism} $(A,m)\to (B,m)$ of $A_\infty$-algebras [GJ] is a
graded linear map $f:A\to B$ commuting with defining cochains of $A$ and
$B$, that is
\[ m_n(f(a_1),\dots , f(a_n)) = f(m_n(a_1,\dots ,a_n)) \]
for all $n$ and all $a_i\in A$.  Note that while associative algebras form a
full subcategory of the cateogry of $A_\infty$-algebras, the inclusion
of DG algebras into the category of $A_\infty$-algebras is not full.

\section{The $X$-machine}
In this section we derive a bicomplex for $A_\infty$-algebras from which
every other (bi)complex to calculate various kinds of Hochschild and
cyclic homology theories for $A_\infty$-algebras can be defined.
This bicomplex is the exact analogue of Connes-Tsygan bicomplex [C$_1$,T],
originally defined for associative algebras.  The key idea here is to extend
Quillen's approach [Q] for associative algebras, to $A_\infty$-algebras.
Once this bicomplex is defined, the rest of the homological algebra of
$A_\infty$-algebras ``follows the book''.  In particular, the
$(b,B)$-bicomplex, the $S$-coperation and Connes' long exact sequence
follow the same pattern as in cyclic homology [C$_1$,C$_2$].

Our main tool to define the cyclic homology of an $A_\infty$-algebra
is the $X$-complex.  The $X$-complex of a (DG)(co) algebra $C$ is
only a first approximation to its various cyclic homology theories.
However, once $X$ is applied to certain universal functors on $C$ one
obtains complexes which are (quasi-) isomorphic to the standard complexes.
One advantage of this approach is that these functors, like the bar
construction, are defined for algebraic structures which are more
flexible than (DG)(co) algebras.

We follow [Q] to define the $X$-complex of a DG coalgebra.  Let $C$ be
a counital coalgebra over a field $K$ of characteristic zero.  Let
$\Omega^1C = \Omega^1_0C$ denote the $C$-bicomodule of universal differential
forms on $C$.  Let
\[ \Omega^1C_\natural = ker\{\Delta_\ell
-\sigma\Delta_r:\Omega^1C\to C\otimes\Omega^1C\} \]
be the subspace of $C$-cocommutators in $\Omega^1C$.  Let
$\partial_0:\Omega^1C_\natural\to C$ be the restriction of $d:\Omega^1C\to C$
to the cocommutator subspace.  Let $\partial_1:C\to\Omega^1C_\natural$ be the
analogue of Hochschild boundary for coalgebras.  The $X$-complex of $C$,
denoted $X(C)$, is the following $\mathbb{Z}$-graded complex which is
2-periodic:
\[ \dots\overset{\partial_1}{\longrightarrow} C
\overset{\partial_0}{\longrightarrow} \Omega^1C_\natural
\overset{\partial_1}{\longrightarrow}
C\overset{\partial_0}{\longrightarrow}\dots\; . \]

Next let $(C,b^\prime )$ be a counital DG coalgebra with \d |b^\prime | =
-1$ and $C = \bigoplus_{i\geq 0}C_i$.  Then one can repeat the above construction
of the $X$-complex to define $X(C)$.  This is a complex in the category
of complexes, i.e. it is a bicomplex in the usual sense.
In fact, more generally, if $\mathscr{C}$ is an abelian tensor cateogry and
$C$ is an algebra or coalgebra object in $\mathscr{C}$, then $X(C)$ is
defined as a complex in $\mathscr{C}$.

For a DG (co) algebra $A$, we define three homologies $XH_\bullet (A)$,
$XC_\bullet (A)$ and $XP_\bullet (A)$ by $XH_\bullet(A) =
H_\bullet (\Omega^1A_\natural )$, $XC_\bullet (A) = H_\bullet
(\mbox{\em Tot}\; X^+(A))$ and $XP_\bullet (A) = H_\bullet (\mbox{\em Tot}\;
X(A))$.

Because of the 2-periodicity in $X(A)$, we obtain a degree 2 map
$S: XC_\bullet (A)\to XC_{\bullet -2}(A)$ and a long exact
sequence
\[ \to XC_n(A)\to XC_{n-2}(A)\to XH_{n-2}(A)\to XC_{n-1}(A)\to\; , \]
similar to Connes' long exact sequence.

A morphism $f:A\to B$ of DG (co) algebras is called an {\em equivalence}
if the induced map $H_\bullet (A,d)\to H_\bullet (B,d)$ is an
isomorphism.

\begin{proposition}
An equivalence $f:A\to B$ induces isomorphisms on $XH_\bullet$,
$XC_\bullet$ and $XP_\bullet$.
\end{proposition}

In particular we can apply this construction to the DG coalgebra $C =
(T^cA[1], b^\prime_m)$, where $(A,m)$ is an $A_\infty$-algebra.  The
corresponding double complex, when $A$ is an associative algebra, is
identified in [Q] and shown to be isomorphic to Connes-Tsygan
bicomplex.  In general, one obtains nothing new for the horizontal
differentials $\partial_0$ and $\partial_1$ except extra signs since we
are working with graded objects.  They are given by $\partial_0 =
1-\lambda$ and $\partial_1 =N=1+\lambda +\dots +\lambda^n$, where
$\lambda :A^{\otimes n+1}\longrightarrow A^{\otimes n+1}$ is
the cyclic shift
\[ \lambda (a_0,\dots ,a_n) = (-1)^\varepsilon(a_n,a_0,\dots ,a_{n-1})\; , \]
where \d \varepsilon = (|a_n|+1)\left(\sum^{n-1}_{i=0}|a_i|+n\right)\;
.$

Let the induced operator on $\Omega^1T^cA[1]_\natural$ be denoted by
$B_m$.  One obtains the following formula for $b_m$:
\begin{align*}
& b_m(a_0,\dots ,a_n) = \sum^{n+1}_{i=0}b_{m_i} (a_0,\dots , a_n)\; ,\\
\intertext{where} & b_{m_i}(a_0,\dots ,a_n) =
b^\prime_{m_i}(a_0,\dots ,a_n) + \sum^n_{j=n-i+1}(-1)^\varepsilon
(m_i(a_j,\dots ,a_{i+j-n-2}),\dots ,a_{j-1})\; .
\end{align*}

Let us denote this double complex by $CC(A,m)$, the part
that is in the first quadrant by $CC^+(A,m)$ and the
zeroth column by $C(A,m)$.  We define the {\em Hochschild
homology} of the $A_\infty$-algebra $(A,m)$ (with coefficients in $A$)
as the homology of $C(A,m)$ and denote it by $HH_\bullet
(A,m)$.  The {\em cyclic homology} of $(A,m)$ is defined as the
homology of the total complex $\mbox{\em Tot}\; CC(A,m)$ and
will be denoted by $HC_\bullet (A,m)$.  We define the {\em periodic
cyclic homology} of $(A,m)$, denoted $HP_\bullet (A,m)$, as the homology
of $\mbox{\em Tot}\; CC(A,m)$, where $\;\hat{}\;$ means we
take direct product instead of direct sum in the total complex.  (The
corresponding homology with direct sums is trivial in all degrees.)

Due to its periodicity, the complex $\mbox{\em Tot}\; CC(A,m)$
has a degree 2 morphism which induces a map $S:HC_\bullet
(A,m)\to HC_{\bullet -2}(A,m)$.  One has the analogue of Connes' long
exact sequence
\[ \longrightarrow HH_\bullet (A,m)\overset{I}{\longrightarrow}
HC_\bullet (A,m)\overset{S}{\longrightarrow}HC_{\bullet -2}(A,m)
\overset{B}{\longrightarrow}HH_{\bullet -1}(A,m)\longrightarrow\; . \]

Let $C^\lambda_n(A,m)=Coker\{ 1-\lambda :C_n(A)\to C_n(A)\}$.  From the
bicomplex relation $b_m(1-\lambda )=(1-\lambda )b^\prime_m$ it is clear
that $(C^\lambda_\bullet (A,m), b_m)$ is a complex.  Moreover, the natural
projection
\[ \mbox{\em Tot}\; CC^+(A,m)\to C^\lambda (A,m) \]
is a quasi-isomorphism.

Defining cyclic cohomology of $A_\infty$-algebras is straightforward.
One should simply dualize the bicomplex $CC(A,m)$ by replacing
tensor products by multilinear functionals.  Let us identify $HC^0
(A,m)$.  A cocycle in $HC^0 (A,m)$ is defined by a {\em closed
graded trace}, that is a linear map $f:A\to k$ such that
$f(ab-(-1)^{|a|\; |b|}ba) = 0$ and $f(m_1a) = 0$ for all $a,b$ in $A$.
Thus $HC^0(A,m)$ is isomorphic to the space of closed graded traces on
$(A,m)$.

Next we turn to the analogue of Connes' operator $B$ and in particular
a $(b,B)$ bicomplex for $A_\infty$-algebras.  This is already achieved
in [GJ] and the relations
\[ B^2 = b_mB+Bb_m=0 \]
verified.  In our approach this comes about as follows.  Let $(A,m)$ be
a {\em unital} $A_\infty$-algebra.  This means there exist an element
$1\in A_0$ such that $m_2(a,1)=m_2(1,a)=a$ for all $a\in A$ and $m_n(a_1,\dots
,a_{i-1},1,a_{i+1}\dots a_n) = 0$ for $n\not= 2$ and all $a_i\in A$.
Let $s:A^{\otimes n}\to A^{\otimes n+1}$ be the standard map
$s(a_1,\dots ,a_n)=(1,a_1,\dots ,a_n)$.  Let $B =
(1-\lambda^{-1})sN:A^{\otimes n}\to A^{\otimes n+1}$.  The
relations $B^2=b_mB+Bb_m=0$ are consequences of bicomplex relations
$N(1-\lambda )=(1-\lambda )N=0$, $b_m(1-\lambda)=(1-\lambda )b^\prime_m$,
$b^2_m=b^{\prime 2}_m = 0$ and the relation
$sb^\prime_m+b^\prime_ms=0$.

The well-known homotopy equivalence between the cyclic and
$(b,B)$-bicomplex carries over to the $A_\infty$-case verbatim.  One can
also consider a normalized $(b,B)$ bicomplex.

We need the following concept and the next proposition for the proof of
theorem 4.4 in Sect. 4.  A morphism $(A,m)\to (B,m)$ of
$A_\infty$-algebras is said to be an {\em equivalence} if the
corresponding map $(T^cA[1],b^\prime_m)\to (T^cB[1],b^\prime_m)$ is a
quasi-isomorphism, i.e. induces an isomorphism $H_\bullet
(T^cA[1],b^\prime_m)\overset{\sim}{\longrightarrow}H_\bullet
(T^cB[1],b^\prime_m)$.

\begin{lemma}
A strict morphism $(A,m)\to (B,m)$ of $A_\infty$-algebras is an
equivalence iff the induced map $H_\bullet (A,m_1)\to H_\bullet (B,m_1)$
is an isomorphism.
\end{lemma}

\begin{proposition}
An equivalence $(A,m)\to (B,m)$ of $A_\infty$-algebras induces
isomorphisms of $HH_\bullet$, $HC_\bullet$, and $HP_\bullet$.
\end{proposition}

\begin{proof}
This is a special case of prop. 2.1 for the DG coalgebra
$T^cA[1]$ and $T^cB[1]$.
\end{proof}

\section{Deformation theory of $A_\infty$-algebras}
The link between deformations of an (associative, commutative, Lie,
etc.) algebra $A$ and the (Hochschild, Harrison, Chevally-Eilenberg,
etc.) cohomology of $A$ with coefficients in $A$ is well known.  One
knows that obstructions for extending a deformation in each order live
in $H^3(A,A)$ and isomorphism classes of deformations are classified by
$H^2(A,A)$.  Moreover, in all of the above cases the Hochschild
cohomology is a {\em Gerstenhaber algebra}, i.e., a graded Poisson algebra.

The link between cyclic cohomology of an associative algebra $A$ and its
deformation theory was first elucidated in [CFS], where it is shown that
if we restrict to deformations that preserve a trace (closed
deformations), then the corresponding obstructions are in $HC^3(A)$ and
$HC^2(A)$.

In this section we define and study the deformation complex of an
$A_\infty$-algebra much in the spirit of the rest of this paper.  We
also establish the link between cyclic cohomology of $A$ and its
deformation complex.  Apparently this fact is of importance in the
cohomology of graph complexes [K,PS].

Let $(A,m)$ be an $A_\infty$-algebra.  We can take
\[ C^\bullet (A,A) = Hom(T^cA[1],A[1])\simeq Coder(T^cA[1],T^cA[1]) \]
as the underlying graded vector space to define the Hochschild cohomology
$H^\bullet (A,A)$.  Note that this is $\mathbb{Z}$-graded, although $A$
is only positively graded.  Also a Lie bracket is defined on $C^\bullet
(A,A)$.  Define a differential
\[ \delta :C^\bullet (A,A)\to C^{\bullet -1}(A,A) \]
by $\delta x = [x,b^\prime_m]$, where we are interpreting $x$ as a
coderivation.  From $b^{\prime 2}_m = 0$ and the Jacobi identity, it
easily follows that $\delta^2=0$ and that $\delta$ is a graded
derivation:
\[ \delta [x,y] = [\delta x,y] + (-1)^{|x|}[x,\delta y]\; , \]
i.e. $(C^\bullet (A,A),\delta ,[\; ,\; ])$ is a differential
$\mathbb{Z}$-graded Lie algebra.  The formulas for the differential and
brackets are as follows.  Let \d f = \sum^\infty_{i=0}f_i$ be a degree
$k$ cochain.  We have
\begin{align*}
\delta f & = \sum^\infty_{n=1}\sum_{i+j=n}(f_i\circ m_j-(-1)m_j\circ f_i)\\
[f_1,f_2] & = \sum^\infty_{n=0}\sum_{i+j=n} (f_i\circ f_j-(-1)f_j\circ
f_i)\; .
\end{align*}

There is, however, more structure hidden in $C^\bullet (A,A)$.  We need
the following simple lemma.

\begin{lemma}
Let $C$ be a graded coalgebra and $(A,m)$ an $A_\infty$-algebra.  Then
there is a natural $A_\infty$-algebra structure on $Hom(C,A)$.
\end{lemma}

\begin{proof}
Use ``multiplications'' on $A$ and comultiplication on $C$ to define
cochains
\[ \widetilde{m}_n:Hom(C,A)^{\otimes n}\longrightarrow Hom(C,A),\quad
n\geq 1\; . \]
Let $\Delta^n:C\to C^{\otimes n}$ be the n-th iteration of the
comultiplication $\Delta :C\to C\otimes C$ of $C$.  Let
$\widetilde{m}_n(f_1,\dots , f_n) = m_n\circ (f_1\otimes\dots\otimes
f_n)\circ\Delta^{(n)}$.  Checking the $A_\infty$ condition is
straightforward.
\end{proof}

The relation between the ``cup product'' and the Lie bracket, even in the
case of associative algebras is quite subtle.  Nevertheless, in [G] it
is shown that $H^\bullet (A,A)$, for $A$ an associative algebra, is a
Gerstenhaber algebra.  That is, the induced cup product in $H^\bullet$
is graded commutative and is compatible with the induced Lie bracket in
the sense that for any $x\in H^\bullet$, the operator $a\mapsto [a,x]$
is a graded derivation of the cup product.

To prove a similar result for $A_\infty$-algebras, we need the notion of
{\em homotopy Gerstenhaber algebra} or $G_\infty$-algebra due to
Gerstenhaber and Voronov [GV].  Let $(B,m)$ be an $A_\infty$-algebra.
A $G_\infty$-structure on $(B,m)$ is an associative product
$T^cB[1]\times T^cB[1]\to T^cB[1]$ on the bar construction of $A$ such
that the codifferential $b^\prime_m$ is a graded derivation of this cup
product.

\begin{lemma}
Let $(B,m)$ be a $G_\infty$-algebra.  Then $H_\bullet (B,m_1)$ is a
Gerstenhaber algebra.
\end{lemma}

\begin{corollary}
Let $(A,m)$ be an $A_\infty$-algebra.  Then the Hochschild cohomology
$H^\bullet (A,A)$ is a Gerstenhaber algebra.
\end{corollary}

At this stage we notice that it is straightforward to define formal
deformations of $A_\infty$-algebras and link it with $H^2(A,A)$ and
$H^3(A,A)$.  Instead we link $HC^\bullet (A)$ to deformations that
preserve a trace, or, equivalently, an invariant bilinear form.
(We are assuming $A$ is unital.)

To illustrate, let us first consider the case where $A$ is an
associative algebra.  We then have the pairings
\begin{align*}
H^p(A,A)\otimes HH^q(A) & \longrightarrow HH^{p+q}(A)\qquad\\
H^k(A,A)\otimes HC^0(A) & \longrightarrow HC^n(A)
\end{align*}
The first map is induced by a morphism of complexes of degree zero
\[ C^\bullet (A,A)\otimes C^\bullet (A)\longrightarrow C^\bullet (A)\; , \]
defined as follows.  For $\varphi :A^{\otimes p}\longrightarrow
A$ in $C^p(A,A)$ and $\tau :A^{\otimes q+1}\longrightarrow K$ in
$C^q(A)$, define $\widetilde{\varphi}:A^{\otimes p+q+1}\longrightarrow K$
by
\[ \tilde{\varphi}(a_0,a_1,\dots ,a_{p+q}) = \tau (a_0\varphi
(a_1,\dots , a_p), a_{p+1},\dots ,a_{p+q})\; . \]
It is straightforward to check that
$\varphi\times\tau\mapsto\widetilde{\varphi}$ is a morphism of complexes and
hence the first pairing is defined.

To define the second pairing, we can interpret $HC^\bullet (A)$ as the
cohomology of $\mbox{\em Tot}\; CC^+(A)$ and define a morphism of complexes
of degree zero
\[ C^\bullet (A,A)\otimes HC^0(A)\longrightarrow \mbox{\em Tot}\; CC^+(A) \]
as follows.  For $\tau :A\to K$ a trace on $A$ and $\varphi :A^{\otimes
n}\longrightarrow A$ in $C^p(A,A)$ define in $\mbox{\em Tot}\; CC^+(A)$ by
\begin{align*}
& \psi = (\psi_p,\psi_{p-1},\dots )\\
& \psi_p(a_0,a_1,\dots ,a_p) = \tau (_0\varphi (a_1,\dots a_p)),\\
& \psi_{k-1}(a_0,\dots ,a_{k-1}) = \tau (\varphi (a_0,\dots ,a_{k-1}))\\
& \psi_i = 0,\qquad\mbox{for}\qquad i<k-1\; .
\end{align*}
Then it is not difficult to check that the above map is a morphism of
complexes.  This means $b\widetilde{\varphi} = \widetilde{\delta\varphi}$
which has already appeared in the first pairing and
$(1-\lambda)\psi_p-b^\prime\psi_{p-1} = (\delta\varphi )_p$,
and this is easy to verify.
This construction has an obvious extension to $A_\infty$ algebras.
Using the second pairing we can transfer cohomological relations in $H^\bullet
(A,A)$ to ones in $HC^\bullet (A)$.  The point is that of course
$HC^\bullet$ is, in general, a ``smaller'' group than $H^\bullet$.  Closed
deformations is an instant where this map can be used.

Let $A$ be an associative algebra and $\tau :A\to K$ a trace on $A$.  A
formal deformation $(A[[t]],\ast )$ defined by
\[ a\ast b = \sum^\infty_{i=0}m_i(a,b)t^i\; , \]
with $m_0(a,b) = ab$, is called {\em closed} (with respect to $\tau$
[CFS]) if
the functional
\begin{align*}
& \widehat{\tau}:A[[t]]\to K[[t]],\\
& \widehat{\tau}\left(\sum^\infty_{i=0}a_it^i\right) =
\sum^\infty_{i=0}\tau(a_i)t^i
\end{align*}
is a trace on $(A[[t]],\ast )$.  It is easy to see that this is
equivalent to
\[ \tau (m_i(a,b)) = \tau (m_i(b,a)) \]
for all $a,b\in A$ and $i\geq 0$.

The question of extending a closed deformation of order $n$ to one of
order $n+1$ amounts to solving the equation
\[ \delta m_{n+1} = m_0\circ m_{n+1}+m_{n+1}\circ m_0 =
-\sum_{i+j=n}m_i\circ m_j\; , \]
for $m_{n+1}$ such that $\tau m_{n+1}$ be symmetric, given that $\tau
m_i$, $0\leq i\leq n$, are symmetric.  This means the right hand side
should represent zero in $H^3(A,A)$.  Using the above pairing, this is
transferred to $HC^3(A)$.  There is a similar argument for equivalence
of closed $\ast$-products.

\begin{proposition}
Let $(A,m)$ be an $A_\infty$-algebra.  Let $\tau\in HC^0(A,m)$ be a
``trace'' on $(A,m)$.  Then the obstructions to extending a closed
$\ast$-deformations of $(A,m)$ within the category of
$A_\infty$-algebras, at any order lie in $HC^3(A,m)$.  Similarly,
extending an equivalence with closed $\ast$-products at any order lies
in $HC^2(A,m)$.
\end{proposition}

\section{Derivations, homotopy invariance and a Goodwillie type theorem}
To prove our main theorem (Theorem 4.5 below), we first extend the
language of derivations and the corresponding Cartan homotopy formula to
$A_\infty$-algebras.  Note that Cartan homotopy formula is the
infinitesimal form of homotopy invariance, from which homotopy
invariance and other results, like Goodwillie's theorem on nilpotent
extensions follows.

Let $(A,m)$ be an $A_\infty$-algebra.  By a {\em graded derivation of
degree} $|D|\in\mathbb{Z}$ of $(A,m)$, we mean a degree $|D|$ map
\[ D:T^cA[1]\to A[1] \]
such that\[ \delta (D) = [m,D] = m\circ D-(-1)^{|D|}\; D\circ m=0\; . \]
Equivalently, $D$ is a cocycle of dimension $|D|$ in the deformation
complex of $A$, introducted in Section 3.  Explicitly, $D$ is defined by
a sequence of maps $D_n:A^{\otimes n}\longrightarrow H$, $n=1,2,\dots$
such that $|D_n| =$\qquad\qquad and
\[ [m,D] = \sum^\infty_{n=2}\sum_n(m_i\circ D_j-(-1)^{|D|}D_j\circ m_i) =
0. \]
In particular $m$ is itself a derivation of degree -1.  It is, however,
trivial as it can be checked that $m=\delta (1)$.

\begin{example}
Let $A$ be a (non-graded) associative algebra considered as an
$A_\infty$-algebra with $A_0=A$, $A_i=\{0\}$, $i\geq 1$, $m_i = 0$ for
$i\not= 2$ and $m_2 =$ multiplication of $A$.  A linear map $D:A^{\otimes
k+1}\longrightarrow A$ defines a map $D:T^cA[1]\to A[1]$ of degree
$-k$ by extending it by zero.  Now $D$ is an $A_\infty$-derivation iff
the original $D$ is a cocycle for the (standard) Hochschild cohomology
$H^{k+1}(A,A)$.  In particular a derivation of algebras in the usual
sense is a derivation in the above sense.
\end{example}

Now a derivation $D:(A,m)\to (A,m)$ of degree $|D|$ induces a
coderivation $b^\prime_D:T^cA[1]\to T^cA[1]$ which is moreover
compatible with the original codifferential $b^\prime_m$:
\[ [b^\prime_D,b^\prime_m] = b^\prime_{[D,m]} = 0\; . \]
So now we have a $DG$ coalgebra $(T^cA[1],b^\prime_m)$ with a
compatible coderivation $b^\prime_D$.

It is well known that derivations of an associative algebra act on
various Hochschild and cyclic complexes of $A$ via the so-called Lie
derivative and one has a Cartan homotopy formula which implies that the
induced action on deRham cohomology and periodic cyclic homology of $A$ is
trivial [C$_1$,Go].  Same formula is also crucial for the proof of
Goodwillie's theorem on periodic cyclic homology of nilpotent
extensions.  In [Go], these results are extended to $DG$ algebras and
applied to singular chains on free loop spaces.

A universal point of view on these matters is as follows.  This will be
very useful in extending these results to $A_\infty$-algebras.  Assume
we have an abelian tensor category $\mathscr{C}$.  One can first extend
the Cuntz-Quillen definition of quasi-free algebras [CQ$_1$,CQ$_2$], to
define a quasi-free algebra or coalgebra object in $\mathscr{C}$.  There are
several equivalent definitions, but the one based on existence of a
connection $\bigtriangledown :\Omega^1A\to\Omega^2A$ for algebra
objects, or a coconnection $\bigtriangledown :\Omega^2C\to\Omega^1C$ for
coalgebra objects is most useful for us.  Here $\Omega^\bullet$ is the
$DG$ algebra of noncommutative differential forms which can be defined
in $\mathscr{C}$.  A connection for algebras is a linear map which is
left $A$-linear and has right Leibniz property $\bigtriangledown (\omega
a)=\bigtriangledown\omega\cdot a+\omega da$.  We call a (co) algebra
object quasi-free if it admits a connection in the above sense.

Now in the above general setting, a derivation $D:C\to C$ induces a {\em
Lie derivative} map $L_D:X(C)\to X(C)$ of degree zero as follows.  In
odd degrees it is simply $D$ itself while in even degrees it is the map
induced by $D$ on $\Omega^1C_\natural$.

The Cartan-homotopy formula of Cuntz and Quillen [CQ$_1$] extends verbatim to
show that if $C$ is quasi-free, then there exists an operator $I_D:X(C)\to X(C)$
of odd degree such that $L_D =\partial I_D+I_D\partial$.

In extending all this to our setup, the category of $DG$ coalgebras, we
have to face the fact that the bar construction is not quasi-free as a
$DG$ coalgebra---if it was, there would be no cyclic homology.  It is,
however, quasi-free, and in fact free, as a graded coalgebra only.  In
[Kh], this problem is addressed and solved by noticing that the vertical
differential in the cyclic bicomplex is given, in even and odd dimension,
by $L_m$, and one obtains a cartan homotopy formula of the form
\[ L_D = [\partial_1+\partial_2,I_D]+I_{\delta D} \]
so that, if $C$ is only quasi-free as a coalgebra and $D$ is compatible
with its differential, then $L_D$ still acts trivially.  The explicit
form of the operator $I_D$ is irrelevant here.

If we specialize the above general result to the case where $C=T^cA[1]$
is the bar construction of an $A_\infty$-algebra $(A,m)$ and $\delta
=b^\prime_D$ is the closed coderivation associated to a derivation $D$
of $(A,m)$, we obtain

\begin{corollary}
The induced map $L_D$ on $HP_\bullet(A,m)$ is zero.
\end{corollary}

\begin{remark}
Similar to $(DG)$ algebras, it is trivial to see that $L_D$ induces the
zero map on $H^{dR}(A,m)$ and also $L_D\circ S$ is zero on $HC_\bullet
(A,m)$.
\end{remark}

Let $(A,m)$ be an $A_\infty$-algebra.  By an $A_\infty$-{\em ideal} we
mean a graded subspace $I\subset A$ such that for all $n\geq 1$,
$m_n(a_1,\dots ,a_n)\in I$ if for some $i$, $a_i\in I$.  Note that in
particular $m_1I\subset I$.  The graded quotient space \d A/I =
\bigoplus^\infty_{i=0}A_i/I_i$ is an $A_\infty$-algebra in a natural way
and the quotient map $A\to A/I$ is a strict morphism of
$A_\infty$-algebras.

\begin{theorem}
Let $(A,m)$ be an $A_\infty$-algebra over a field of characteristic zero
and let $I\subset A$ be an $A_\infty$-ideal such that $I_0 = 0$.  Then
the quotient map $(A,m)\to (A/I,m)$ induces an isomorphism
\[ HP_\bullet (A,m)\longrightarrow HP_\bullet (A/I,m)\; . \]
\end{theorem}

\begin{proof}
As in [Go], let
\[ gr(A,I,m) = \bigoplus_{k\geq 0}I^k/I^{k+1}\; . \]
This is an $A_\infty$-algebra.  There is a derivation, the so-called
{\em number operator}, acting on $gr(A,I,m)$ by multiplying a
homogeneous element $a$ by its degree $|a|$.  The cyclic complex of
$gr(A,I,m)$ can be identified as follows.  Filter the cyclic complex
$CC_\bullet (A,m)=\mbox{\em Tot}\; X^+(T^CA[1])$ by subcomplexes $F^k$, $k\geq
0$, where for each $n$, the $n$-tensor components of $F^k$ are tensors
of the form
\[ \bigoplus_{k_0+\cdots +k_n=k}I^{k_0}\otimes\cdots\otimes I^{k_n}\; .
\]
Then the cyclic complex of $gr(A,I,m)$ is given by
\[ CC_\bullet (gr(A,I,m))\simeq\bigoplus_{k\geq 0}F^k/F^{k+1}\; . \]
The operator $S$ descends to subcomplexes $F^k$.  It suffices to show
that the map
\[ H_{n+2k}(F^\prime )\overset{S^k}{\longrightarrow}H_n(F^\prime ) \]
is zero for $n\leq k$.  As in [Go], this follows from the following
two observations:

\begin{enumerate}
\item[1.] $S^k:H_n(F^\prime/F^{k+1})\longrightarrow H_n(F^\prime
/F^{k+1})$ is zero for all $n$.
\item[2.] $H_n(F^k)=0$ for $n<k$.  This is simply true because $F^k$ has
non chains in dimensions less than $k$.
\end{enumerate}
\end{proof}

Finally, we prove our main result:

\begin{theorem}
Let $(A,m)$ be an $A_\infty$-algebra over a field of characteristic
zero.  Let $H_\bullet (A) = H_\bullet (A,m_1)$ be the homology of $A$ in
degree zero with its induced associative product.  Then
$HC^{per}_\bullet (A,m)\simeq HC^{per}_\bullet (H_0(A))$.  [Note that in
the right hand side we have the ordinary cyclic homology functor.]
\end{theorem}

\begin{proof}
This is a consequence of Lem. 2.2, Prop. 2.3 and Prop. 4.4.  Define a
new $A_\infty$-algebra \d B = \bigoplus_{i\geq 0}B_i$ by $B_0 = A_0/Im\;
m_1$ and $B_i = A_i$, $i\geq 1$.  Note that $B_0 = H_0(A,m_1)$.  We have
a strict morphism of $A_\infty$-algebras $(A,m)\to (B,m)$.  By Lem. 2.2,
this is an equivalence and hence induces an isomorphism on $HP_\bullet$.
Let \d I = \bigoplus_{i\geq 1}B_i$.  Then $(B,I)$ satisfies conditions of
Prop. 2.3 and hence we have an isomorphism $HP_\bullet (B,m)\to HP_\bullet
(B_0)$.  It follows that the map $HP_\bullet (A,m)\to HP_\bullet
(H_0(A,m_1))$, being a composition of two isomorphisms, is an
isomorphism.
\end{proof}

\begin{remark}
Although we have not checked it, but one can perhaps prove a stronger
result.  Let us call a morphism of $A_\infty$-algebras $(A,m)\to (B,m)$
to be {\em 1-connected} if it induces isomorphism on $H_0$ and is
surjective on $H_1$.  Then a 1-connected map of $(A,m)\to
(B,m)$ of $A_\infty$-algebras induces an isomorphism
\[ HP_\bullet (A,m)\overset{\sim}{\longrightarrow}
HP_\bullet (B,m)\; . \]
If $A$ and $B$ are $DG$ algebras, this is Theorem IV.2.1 of [Go].
\end{remark}

\vskip30pt

\noindent Masoud Khalkhali\\
University of Western Ontario\\
London, Canada\\
N6A 5B7\\
masoud@julian.uwo.ca

\end{document}